# Asymptotics: Particles, Processes and Inverse Problems

**Festschrift for Piet Groeneboom**

**Eric A. Cator, Geurt Jongbloed, Cor Kraaikamp,**
**Hendrik P. Lopuhaä, Jon A. Wellner, Editors**

# Contents

# Preface

In September 2006, Piet Groeneboom officially retired as professor of statistics at Delft University of Technology and the Vrije Universiteit in Amsterdam. He did so by delivering his farewell lecture 'Summa Cogitatio' ([42] in Piet's publication list) in the Aula of the university in Delft. To celebrate Piet's impressive contributions to statistics and probability, the workshop 'Asymptotics: particles, processes and inverse problems' was held from July 10 until July 14, 2006, at the Lorentz Center in Leiden. Many leading researchers in the fields of probability and statistics gave talks at this workshop, and it became a memorable event for all who attended, including the organizers and Piet himself.

This volume serves as a Festschrift for Piet Groeneboom. It contains papers that were presented at the workshop as well as some other contributions, and it represents the state of the art in the areas in statistics and probability where Piet has been (and still is) most active. Furthermore, a short CV of Piet Groeneboom and a list of his publications are included.

Eric Cator
Geurt Jongbloed
Cor Kraaikamp
Rik Lopuhaä
Delft Institute of Applied Mathematics
Faculty of Electrical Engineering,
Mathematics and Computer Science
Delft University of Technology
The Netherlands

Jon Wellner
Department of Statistics
University of Washington, Seattle
USA

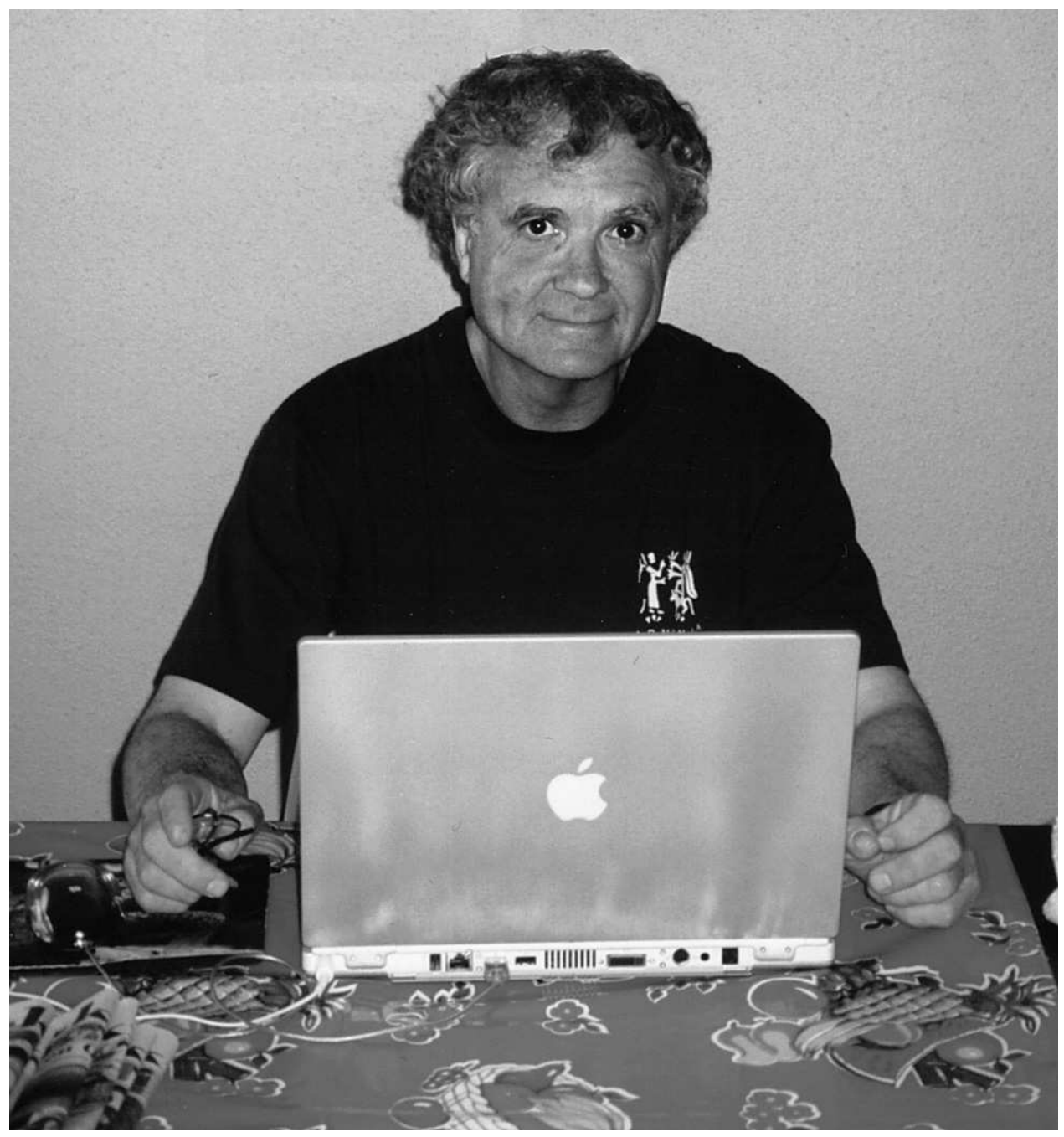

Piet in a characteristic pose. Amsterdam, 2003.

# Curriculum Vitae of Piet Groeneboom

| | |
|---|---|
| **Born:** | September 24, 1941, The Hague |
| **Citizenship:** | The Netherlands |
| **Dissertation:** | *Large deviations and asymptotic efficiencies*<br>1979, Vrije Universiteit, Amsterdam.<br>supervisor: J. Oosterhoff. |

**Professional Career:**

| | | |
|---|---|---|
| Mathematical Centre (MC) | Researcher and consultant | 1973–1984 |
| University of Washington, Seattle | Visiting assistant professor | 1979–1981 |
| University of Amsterdam | Professor of statistics | 1984–1988 |
| Delft University of Technology | Professor of statistics | 1988–2006 |
| Stanford University | Visiting professor | 1990 |
| Université Paris VI | Visiting professor | 1994 |
| University of Washington, Seattle | Visiting professor | 1998,<br>1999, 2006 |
| University of Washington, Seattle | Affiliate professor | 1999– |
| Vrije Universiteit Amsterdam | Professor of statistics | 2000–2006 |
| Institut Henri Poincaré, Paris | Visiting professor | 2001 |

**Miscellanea:**

Rollo Davidson prize 1985, Cambridge UK.

Fellow of the IMS and elected member of ISI.

Visitor at MSRI, Berkeley, 1983 and 1991.

Three times associate editor of *The Annals of Statistics*.

Invited organizer of a DMV (Deutsche Mathematiker Vereinigung) seminar in Günzburg, Germany, 1990.

Invited lecturer at the Ecole d'Eté de Probabilités de Saint-Flour, 1994.

## Publications of Piet Groeneboom

*April 2007*

1. Rank tests for independence with best strong exact Bahadur slope (with Y. Lepage and F.H. Ruymgaart), *Zeitschrift für Wahrscheinlichkeitstheorie und Verwandte Gebiete* **36** (1976), 119–127.
2. Bahadur efficiency and probabilities of large deviations (with J. Oosterhoff), *Statist. Neerlandica* **31** (1977), 1–24.
3. Relevant variables in the advices of elementary school teachers on further education; an analysis of correlational structures (in Dutch, with J. Hoogstraten, G.J. Mellenbergh and J.P.H. van Santen), *Tijdschrift voor Onderwijsresearch (Journal for Educational Research)* **3** (1978), 262–280.
4. Large deviation theorems for empirical probability measures (with J. Oosterhoff and F.H. Ruymgaart), *Ann. Probability* **7** (1979), 553–586.
5. *Large deviations and asymptotic efficiencies*, Mathematical Centre Tract **118** (1980), Mathematical Centre, Amsterdam
6. Large deviations of goodness of fit statistics and linear combinations of order statistics (with G.R. Shorack), *Ann. Probability* **9** (1981), 971–987.
7. Bahadur efficiency and small-sample efficiency (with J. Oosterhoff), *Int. Statist. Rev.* **49** (1981), 127–141.
8. The concave majorant of Brownian motion, *Ann. Probability* **11** (1983), 1016–1027.
9. Asymptotic normality of statistics based on convex minorants of empirical distribution functions (with R. Pyke), *Ann. Probability* **11** (1983), 328–345.
10. Estimating a monotone density, in *Proceedings of the Conference in honor of Jerzy Neyman and Jack Kiefer, Vol. II* (Eds. L.M. Le Cam and R.A. Olshen), Wadsworth, Inc, Belmont, California (1985), 539–555.
11. Some current developments in density estimation, in *Mathematics and Computer Science, CWI Monograph 1* (Eds. J.W. de Bakker, M. Hazewinkel, J.K. Lenstra), Elsevier, Amsterdam (1986), 163–192.
12. Asymptotics for incomplete censored observations, Mathematical Institute, University of Amsterdam (1987), Report 87-18.
13. Limit theorems for convex hulls, *Probab. Theory Related Fields* **79** (1988), 327–368.
14. Brownian motion with a parabolic drift and Airy functions, *Probab. Theory Related Fields* **81** (1989), 79–109.
15. Discussion on "Age-specific incidence and prevalence, a statistical perspective", by Niels Keiding in the *J. Roy. Statist. Soc. Ser. A.* **154** (1991), 371–412.
16. *Information bounds and nonparametric maximum likelihood estimation* (with J.A. Wellner), Birkhäuser Verlag (1992).
17. Discussion on "Empirical functional and efficient smoothing parameter selection" by P. Hall and I. Johnstone in the *J. Roy. Statist. Soc. Ser. B.* **54** (1992), 475–530.
18. Isotonic estimators of monotone densities and distribution functions: basic facts (with H.P. Lopuhaä), *Statist. Neerlandica* **47** (1993), 175–183.
19. Flow of the Rhine river near Lobith (in Dutch: "Afvoertoppen bij Lobith"), in *Toetsing uitgangspunten rivierdijkversterkingen, Deelrapport 2: Maatgevende belastingen* (1993), Ministerie van Verkeer en Waterstaat.
20. Limit theorems for functionals of convex hulls (with A.J. Cabo), *Probab. Theory Related Fields* **100** (1994), 31–55.

21. Nonparametric estimators for interval censoring, in *Analysis of Censored Data* (Eds. H. L. Koul and J. V. Deshpande), IMS Lecture Notes-Monograph Series **27** (1995), 105–128.
22. Isotonic estimation and rates of convergence in Wicksell's problem (with G. Jongbloed), *Ann. Statist.* **23** (1995), 1518–1542.
23. Computer assisted statistics education at Delft University of Technology, (with de P. Jong, D. Tischenko and B. van Zomeren), *J. Comput. Graph. Statist.* **5** (1996), 386–399.
24. Asymptotically optimal estimation of smooth functionals for interval censoring, part 1 (with R.B. Geskus), *Statist. Neerlandica* **50** (1996), 69–88.
25. Lectures on inverse problems, in *Lectures On Probability and Statistics*. Ecole d'Eté de de Probabilités de Saint-Flour XXIV (Ed. P. Bernard), Lecture Notes in Mathematics **1648** (1996), 67–164. Springer Verlag, Berlin.
26. Asymptotically optimal estimation of smooth functionals for interval censoring, part 2 (with R.B. Geskus), *Statist. Neerlandica* **51** (1997), 201–219.
27. Extreme Value Analysis of North Sea Storm Severity (with C. Elsinghorst, P. Jonathan, L. Smulders and P.H. Taylor), *Journal of Offshore Mechanics and Arctic Engineering* **120** (1998), 177–184.
28. Asymptotically optimal estimation of smooth functionals for interval censoring, case 2 (with R.B. Geskus), *Ann. Statist.* **27** (1999), 627–674.
29. Asymptotic normality of the $L_1$-error of the Grenander estimator (with H.P. Lopuhaä and G. Hooghiemstra), *Ann. Statist.* **27** (1999), 1316–1347.
30. Integrated Brownian motion conditioned to be positive (with G. Jongbloed and J.A. Wellner), *Ann. Probability* **27** (1999), 1283–1303.
31. A monotonicity property of the power function of multivariate tests (with D.R. Truax), *Indag. Math.* **11** (2000), 209–218.
32. Computing Chernoff's distribution (with J.A. Wellner), *J. Comput. Graph. Statist.* **10** (2001), 388–400.
33. A canonical process for estimation of convex functions: the "invelope" of integrated Brownian motion $+\,t^4$ (with G. Jongbloed and J.A. Wellner), *Ann. Statist.* **29** (2001), 1620–1652.
34. Estimation of convex functions: characterizations and asymptotic theory (with G. Jongbloed and J.A. Wellner), *Ann. Statist.* **29** (2001), 1653–1698.
35. Ulam's problem and Hammersley's process, *Ann. Probability* **29** (2001), 683–690.
36. Hydrodynamical methods for analyzing longest increasing subsequences, *J. Comput. Appl. Math.* **142** (2002), 83–105.
37. Kernel-type estimators for the extreme value index (with H.P. Lopuhaä and P.-P. de Wolf), *Ann. Statist.* **31** (2003), 1956–1995.
38. Density estimation in the uniform deconvolution model (with G. Jongbloed), *Statist. Neerlandica* **57** (2003), 136–157.
39. Hammersley's process with sources and sinks (with E.A. Cator), *Ann. Probability* **33** (2005), 879–903.
40. Second class particles and cube root asymptotics for Hammersley's process (with E.A. Cator), *Ann. Probability* **34** (2006), 1273–1295.
41. Estimating the upper support point in deconvolution (with L.P. Aarts and G. Jongbloed). To appear in the *Scandinavian journal of Statistics*, 2007.
42. Summa Cogitatio. To appear in *Nieuw Archief voor Wiskunde* (magazine of the Royal Dutch Mathematical Association) (2007).
43. Convex hulls of uniform samples from a convex polygon, Conditionally accepted for publication in *Probability Theory and Related Fields*.

44. Current status data with competing risks: Consistency and rates of convergence of the MLE (with M.H. Maathuis and J.A. Wellner). To appear in *Ann. Statist.* (2007).
45. Current status data with competing risks: Limiting distribution of the MLE (with M.H. Maathuis and J.A. Wellner). To appear in *Ann. Statist.* (2007).
46. The support reduction algorithm for computing nonparametric function estimates in mixture models (with G. Jongbloed and J.A. Wellner). Submitted.

# Contributors to this volume

Balabdaoui, F., *Université Paris-Dauphine*
Birgé, L., *Université Paris VI*
Boysen, L., *Universität Göttingen*

Dümbgen, L., *University of Bern*

Ferrari, P. A., *Universidade de São Paulo*

van de Geer, S. A., *ETH Zürich*
Gill, R. D., *Leiden University*

Huang, J., *University of Iowa*

Lalley, S. P., *University of Chicago*
Liebscher, V., *Universität Greifswald*

Mammen, E., *Universität Mannheim*
Munk, A., *Universität Göttingen*

Pollard, D., *Yale University*

Rufibach, K., *University of Bern*

Seppäläinen, T., *University of Wisconsin-Madison*
Sisko, V. V., *Universidade de São Paulo*

van der Vaart, A. W., *Vrije Universiteit Amsterdam*

Wellner, J. A., *University of Washington*
Wittich, O., *Technische Universiteit Eindhoven*

Xie, H., *University of Iowa*

Yu, K., *Universität Mannheim*